\documentclass[12pt]{article}
\usepackage{amsmath}
\usepackage{amssymb}
\usepackage{a4}
\usepackage{theorem}
\usepackage{amscd}
\usepackage[mathscr]{eucal}
\input{diagrams}
\newcommand{\nx}{\negthinspace}

\newcommand{\oA}{\overline{A}}
\newcommand{\oB}{\overline{B}}
\newcommand{\tB}{\tilde{B}}
\newcommand{\CC}{\mathbb{C}}
\newcommand{\oC}{\overline{C}}
\newcommand{\oCB}{\oC_B}

\newcommand{\oD}{\overline{D}}
\newcommand{\tD}{\tilde{D}}
\newcommand{\FF}{\mathbb{F}_2}
\newcommand{\tG}{\tilde{G}}
\newcommand{\tH}{\tilde{H}}
\newcommand{\oI}{\overline{I}}
\newcommand{\oJ}{\overline{J}}
\newcommand{\sL}{\mathscr{L}}
\newcommand{\sM}{\mathscr{M}}

\newcommand{\PP}{\mathbb{P}}
\newcommand{\sP}{\mathscr{P}}

\newcommand{\tPP}{\tilde{\mathbb{P}}}
\newcommand{\QQ}{\mathbb{Q}}

\newcommand{\oS}{\overline{S}}
\newcommand{\oSd}{\smash{\overline{S}}'}
\newcommand{\oSdd}{\smash{\overline{S}}''}
\newcommand{\sS}{\mathscr{S}}
\newcommand{\tS}{\tilde{S}}
\newcommand{\oT}{\overline{T}}

\newcommand{\hV}{\hat{V}}

\newcommand{\tV}{\tilde{V}}
\newcommand{\oX}{\overline{X}}
\newcommand{\tX}{\tilde{X}}
\newcommand{\ZZ}{\mathbb{Z}}
\newcommand{\tpi}{\tilde{\pi}}

\newcommand{\tsigma}{\tilde{\sigma}}
\newcommand{\ttau}{\tilde{\tau}}
\newcommand{\ww}{\left|w\right|}

\newcommand{\Aut}{\operatorname{Aut}}
\newcommand{\codim}{\operatorname{codim}}

\newcommand{\dimF}{\operatorname{dim}_{\FF}}

\newcommand{\Pic}{\operatorname{Pic}}

\newcommand{\minors}{\operatorname{minors}}
\newcommand{\mult}{\operatorname{mult}}

\newcommand{\im}{\operatorname{im}}

\newcommand{\sing}{\operatorname{sing}}

\newcommand{\coker}{\operatorname{coker}}
\newcommand{\cohom}[2]{H^{#1}\!\left(#2\right)}

\newcommand{\cohomd}[2]{h^{#1}\!\left(#2\right)}

\newcommand{\cohomtPP}{\cohom{2}{\smash{\tPP}_3,\ZZ}}
\newcommand{\cohomtV}{\cohom{2}{\smash{\tV},\ZZ}}
\newcommand{\cohomstV}{\cohom{\ast}{\smash{\tV},\ZZ}}
\newcommand{\cohomtB}{\cohom{2}{\smash{\tB},\ZZ}}
\newcommand{\cohomtX}{\cohom{2}{\smash{\tX},\ZZ}}
\newcommand{\cohomdtX}{\cohomd{2}{\smash{\tX},\ZZ}}
{\theorembodyfont{\rm}
    \newtheorem{definition}{Definition}[section]
    \newtheorem{remark}[definition]{Remark}
    \newtheorem{example}[definition]{Example}
}
\newtheorem{lemma}[definition]{Lemma}

\newtheorem{proposition}[definition]{Proposition}
\newtheorem{theorem}[definition]{Theorem}
\newcommand{\proof}{{\noindent\bf Proof:\ }}

\newcommand{\proofend}{$\square$\medskip\par}
\newcommand{\bibauthor}[2]{{\it{#2}~{#1},}}
\newcommand{\bibtitlea}[1]{#1,}            
\newcommand{\bibtitleb}[1]{{\it #1},}      
\newcommand{\bibcompany}[1]{#1}
\newcommand{\bibtown}[1]{#1}
\newcommand{\bibthesis}{Ph.D.\ tthesis}
\newcommand{\bibyear}[1]{(#1),}
\newcommand{\bibyearx}[1]{(#1)}
\newcommand{\bibjournal}[2]{#1\ {\bf #2}}
\newcommand{\bibpages}[2]{{#1}--{#2}}
\newcommand{\bibend}{.}

\begin{document}
\newcounter{hlp}
\ \vspace*{5cm}\\
\noindent{\LARGE\bf On the divisor class group of double solids\\[2ex]}
{\large\bf Stephan Endra\ss\\[1ex]}
Johannes Gutenberg-Universit\"at, Mainz, Germany\\[2ex]
\today\\[5ex]
{\bf Abstract}\\[1ex]
{\small For a double solid $V\rightarrow \PP_3$ branched
over a surface $B\subset\PP_3\left(\CC\right)$ with only ordinary
nodes as singularities, we give a set of generators
of the divisor class group
$\Pic\left(\smash{\tV}\right)\cong\cohomtV$ in terms
of contact surfaces of $B$ with only superisolated
singularities in the nodes of $B$.
As an application we give a condition when
$\cohomstV$ has no 2-torsion. All possible cases
are listed if $B$ is a quartic.
Furthermore we give a new lower bound for
the dimension of the code of $B$.
\setcounter{section}{-1}
\section{Introduction}\label{sect:intro}
In this note we consider threefolds $V$ given by
$\pi\colon V\rightarrow \PP_3\left(\CC\right)$, where $\pi$ is
the  double cover branched exactly over a nodal surface 
$B\subset\PP_3$ of even degree $b$.
A nodal surface is a surface which is smooth except finitely
many  ordinary double points (nodes).
These double covers $V$ of $\PP_3$ are called {\em double solids}.
$V$ has ordinary double points over the nodes of $B$. They
can be resolved by a small resolution $\hV$ or a big
resolution $\tV$. The threefolds $V$, $\hV$, $\tV$ and their relations
have been studied in detail in \cite{clemens},\cite{werner}. There
is a commutative diagram
\begin{diagram}[height=2.5em,width=2.5em]
        \tV & \rTo^{\tpi} & \tPP_3 & \lInto & \tB \\
        \dTo^{\tsigma} & & \dTo^{\sigma}& & \dTo^{\left.\sigma\right|_{\tB}}\\
        V & \rTo^{\pi}& \PP_3 & \lInto &  B
\end{diagram}
where $\sigma\colon\tPP_3\rightarrow\PP_3$ is the blowup of
$\PP_3$ in the nodes of $B$. For every variety $X\subset\PP_3$,
denote by $\tX$ the proper transform of $X$ with respect to $\sigma$.
Then $\tpi$ is the double
cover of $\tPP_3$ branched exactly over $\tB$ and
$\tsigma$ is the big resolution of all double points of $V$.
The double cover $\pi$ (resp.~$\tpi$) defines an involution
$\tau$ on $V$ (resp.~$\ttau$ on $\tV$) exchanging sheets.
Let $H$ be the class of a general hyperplane in $\PP_3$
and denote by $\sS$ the singular locus of $B$.
Let $\mu=\left|\sS\right|$.
For any node $P\in\sS$, let $E_P=\sigma^{-1}(P)\cong\PP_2$ the
corresponding exceptional divisor.
For any set of nodes $N\subseteq\sS$, let
$E_N=\sum_{P\in N}E_P$.
Then $\cohomtPP$ is generated freely by the classes of $\tH$ and
$E_P$, $P\in\sS$. Moreover
$\cohomtV=\ZZ\oplus\ZZ^{\mu}\oplus\ZZ^d$ is free of rank $1+\mu +d$,
where $d$ is called {\em defect of $V$}. The first
two summands are just the image of the induced (injective) map
$\tpi^\ast\colon\cohomtPP\rightarrow\cohomtV$. The third summand is
explained in this note.
The defect $d$ of $V$ can be computed \cite[3.16]{clemens} as
\begin{equation}\label{equation:defect}
    d = \dim\sM - \left(\left( 3n/2-1 \atop 3\right)-\mu\right),
\end{equation}
where $\sM$ is the $\CC$ vector space of all homogeneous polynomials
of degree $3n/2-4$ vanishing on $\sS$.
Note that $d$ is the difference between the actual and the
estimated dimension of $\sM$.
 So a high defect
indicates that the nodes of $B$ are in a very special position
with respect to polynomials of degree $3n/2-4$.

There is another number which measures how special the position
of the nodes of $B$ is.
It is the dimension of the code of $B$.
The code of $B$ is the $\FF$ vector space
\begin{equation*}
    \oC_B = \ker\left(
    \cohom{2}{\smash{\tPP_3},\FF}\overset{rest.}{\longrightarrow}
    \cohom{2}{\smash{\tB},\FF}\right).
\end{equation*}
Projecting $\oC_B$ onto the second sumand of
$\cohom{2}{\smash{\tPP_3},\FF}=\FF\oplus\FF^\mu$ induces an injection
$\left(\oC_B,+\right)\rightarrow\left(\sP(\sS),\Delta\right)$
where $\Delta$ denotes the symmetric difference of
subsets of $\sS$. A subset $N\subseteq\sS$ is called
strictly (weakly) even, if the class
of $E_N$ ($\tH+E_N$) is divisible by 2 in
$\cohomtB$. In the sequel we identify the even sets of
nodes with the elements of $\oC_B$.
The even sets of nodes $w\in\oC_B$ are in
correspondence to surfaces $S\subset\PP_3$ which have
contact to $B$, i.e.~which satisfy $S.B=2D$ for some
curve $D$ on $B$. Then $w$ is just the set of points in $\sS$
where $D$ fails to be cartier on $B$. In this situation
$S$ is called a contact surface of $B$ which cuts out
$w$ via $D$. Note that all surfaces cutting out $w$
on $B$ have even (odd) degree if $w$ is strictly (weakly) even.
For a detailed description of even sets of nodes, we refer to
\cite{catanese} and \cite{endrassnodes}.

In the first section we recall some well known properties of
double covers. In the second section we construct generators
of the third summand $\ZZ^d$ of $\cohomtV$
(theorem \ref{theorem:classes}). This is done by
identifying the classes of contact surfaces of $B$ which
split under the double cover $\tpi$. As an application,
we study double solids branched along nodal quartic surfaces
(theorem \ref{theorem:quartics}).
The case of a 6-nodal quartic
with defect 1 had already been studied in \cite[lemma 3.25]{clemens}.
Finally we give a new lower bound for the dimension of the
code of $B$ (proposition \ref{proposition:dim})
which slightly improves the bound of Beauville
\cite[p.~210]{beauville}.

I would like to thank S.~Cync for helpful discussions on
this topic.
\section{Double covers}
Let $\pi\colon X\rightarrow Y$ be a double cover of smooth
$n$-dimensional varieties branched along a smooth divisor
$D$ on $Y$. For any abelian group $G$, $\pi$ induces an
injection of sheaves
$i_G\colon \tG_Y\rightarrow\pi_{\ast}\tG_X$. Let $Q$
be the $\QQ$-sheaf defined by the exact sequence
\begin{diagram}[height=2em,width=2em]
    0 & \rTo & \QQ_Y & \rTo^{i_\QQ} & \pi_{\ast}\QQ_X & \rTo & Q & \rTo & 0.
\end{diagram}
Averaging over the fibre gives a map
$p_\QQ\colon \pi_{\ast}\QQ_X\rightarrow\QQ_Y$ with $p_\QQ\circ i_\QQ =1$,
hence the above sequence splits. Consequently all induced maps
$\pi^{\ast}\colon\cohom{i}{Y,\QQ}\rightarrow\cohom{i}{X,\QQ}$ are
injective. Hence if $\cohom{\ast}{Y,\ZZ}$ contains no torsion,
then also all maps $\pi^{\ast}\colon\cohom{i}{Y,\ZZ}\rightarrow\cohom{i}{X,\ZZ}$
are injective.

Now consider $G=\FF$. Let $F$ be the sheaf defined by
\begin{equation}\label{seq:1}
\begin{diagram}[height=2em,width=2em]
    0 & \rTo & {\FF}_Y &\rTo^{i_{\FF}} & \pi_\ast{\FF}_X & \rTo & F & \rTo & 0.
\end{diagram}
\end{equation}
This sequence does not split in general, but we have the following
\begin{lemma}
    The $\FF$-sheaf $F$ fits into the exact sequence
    \begin{equation}\label{seq:2}
    \begin{diagram}[height=2em,width=2em]
        0 & \rTo & F & \rTo & {\FF}_Y & \rTo & {\FF}_D & \rTo & 0.
    \end{diagram}
    \end{equation}
\end{lemma}
\proof Let $F'$ be the sheaf defined by the exact sequence (\ref{seq:2}).
For all $y\in Y$ we have
\begin{equation*}
    F'_y=\left\{\begin{array}{c@{\quad}l}
        \FF & \text{if $y\not\in D$,}\\
        0   & \text{if $y\in D$}
    \end{array}\right\}=F_y.
\end{equation*}
So we can find an open cover $\left\{U_\alpha\right\}_{\alpha\in A}$
of the variety $Y$ and a family of isomorphisms $\left\{f_\alpha\colon
F\left(U_\alpha\right)\rightarrow F'\left(U_\alpha\right)\right\}_{\alpha\in A}$.
But $\Aut\left(\FF\right)=1$, hence all $f_\alpha$ glue together.
This implies $F\cong F'$.
\proofend
Let us consider our initial situation $X=\tV$, $Y=\tPP_3$ and $D=\tB$.
From (\ref{seq:1}) and (\ref{seq:2}) we get exact sequences
\begin{diagram}[height=2em,width=2em]
    \cohom{1}{\smash{\tB},\FF} & \rTo &
    \cohom{2}{\smash{\tPP_3},F} & \rTo &
    \cohom{2}{\smash{\tPP_3},\FF} & \rTo &
    \cohom{2}{\smash{\tB},\FF},\\
    \cohom{2}{\smash{\tPP_3},\FF} & \rTo &
    \cohom{2}{\smash{\tV},\FF} & \rTo &
    \cohom{2}{\smash{\tPP_3},F} & \rTo &
    \cohom{3}{\smash{\tPP_3},\FF}.
\end{diagram}
But the integral cohomology of $\tPP_3$ and $\tB$ are torsion free
with $\cohom{1}{\smash{\tB},\ZZ}=\cohom{3}{\smash{\tPP_3},\ZZ}=0$.
So by the universal coefficient theorem
$\cohom{1}{\smash{\tB},\FF}=\cohom{3}{\smash{\tPP_3},\FF}=0$.
This implies that
\begin{align*}
    \oC_B & =\ker\left( \cohom{2}{\smash{\tPP_3},\FF}\overset{rest.}{\longrightarrow}
    \cohom{2}{\smash{\tB},\FF}\right)
    \cong\cohom{2}{\smash{\tPP_3},F} \\
     & \cong
    \coker\left(\cohom{2}{\smash{\tPP_3},\FF}\overset{\tpi^\ast}{\longrightarrow}
    \cohom{2}{\smash{\tV},\FF}\right).
\end{align*}
In general $\cohomstV$ has torsion, so the universal coefficient theorem gives
$\cohom{2}{\smash{\tV},\FF}\cong\left(\cohomtV\otimes_\ZZ \FF\right)\oplus T_2$,
where $T_2$ is the 2-torsion of $\cohom{3}{\smash{\tV},\ZZ}$ (or equivalently
the 2-torsion of
$\cohom{4}{\smash{\tV},\ZZ}$). So we get the following
\begin{lemma}\label{lemma:dim}
    $\oC_B\cong\FF^d\oplus T_2$, where $T_2$ is the 2-torsion of
    $\cohom{3}{\smash{\tV},\ZZ}$. In particular
    $\dim_{\FF}\oC_B\geq d$ with equality iff $\cohomstV$ has no
    2-torsion.
\end{lemma}
This simple lemma enables us in section 3 to give a necessary and
sufficient condition for $\cohomstV$ to have no 2-torsion.
\section{Generators of $\cohomtV$}
Consider the induced map
$\tpi^\ast\colon\cohomtPP\rightarrow\cohomtV$ and let
$C=\coker\tpi^\ast$, $I=\im\tpi^\ast$.
In this section we explicitely give a map $w\colon C\rightarrow\oC_B$
which assigns to a class $\alpha\in C$ a surface $S\subset\PP_3$
which cuts out an even set of nodes $w(\alpha)\in\oC_B$.
The induced map $w_2:(C\otimes_\ZZ\FF)\rightarrow\oC_B$ is seen to be
an injection. Moreover there is a geometric explanation
of the image of $w_2$.

Remember that for any variety $X\subset\PP_3$, its proper transform
with respect to $\sigma$ is denoted by $\tX$. Let 
$\oX=\tpi^{-1}(\smash{\tX})$. Then $\oX$ defines a cohomology
class $[\oX]\in\cohom{\ast}{\smash{\tV},\ZZ}$. 
We need a preparatory lemma:
\begin{lemma}\label{lemma:cut}
    Let $w\in\oC_B$ be cut out by $S$ via $D$. If $\tS.\tB=2\tD$, then
    $w=0\Leftrightarrow [\tS]\in2\cohomtPP$.
\end{lemma}
\proof For all $P\in\sS$ the integer
$\mult(D,P)=\tD.E_P$ is even (odd) if $P\not\in w$ ($P\in w$).
Moreover $\tS\sim_{lin}s\tH-\sum_{P\in\sS}a_PE_P$ for some
integers $s, a_P$, $P\in\sS$. Let $C_P$ be the exceptional conic
on $\tB$ corresponding to $P\in\sS$. Note that $C_P$ has
self intersection $C_P^2=-2$. 
Then
$\left.\smash{\tS}\right|_{\tB}=2\tD\sim_{lin}s\left.\smash{\tH}\right|_{\tB}-
\sum_{P\in\sS}a_PC_P$. It follows that
$\tD.E_P=(1/2)a_PC_P.E_P=-a_P$, so 
\begin{align*}
    w=0 &\Longleftrightarrow 2\mid \tD.E_P\quad\forall P\in\sS
          \quad\text{and}\quad 2\mid s\\
        &\Longleftrightarrow 2\mid a_P 
         \quad\forall P\in\sS\quad\text{and}\quad 2\mid s \\
        &\Longleftrightarrow [\tS]\in2\cohomtPP.
\end{align*}
\proofend
\begin{proposition}\label{proposition:alpha}
    For every $\alpha\in C$, there exists a surface $S\subset\PP_3$
    satisfying
    \begin{itemize}
        \item[A1)] $S.B=2D$ ($S$ and $B$ have contact along a curve $D$),
        \item[A2)] $\tS$ and $\tD$ are smooth with $\tS.\tB=2\tD$,
        \item[A3)] $[\oS]=\beta+\ttau^\ast\beta$ splits into
            conjugate classes with $\beta\equiv\alpha\bmod I$.
    \end{itemize}
\end{proposition}
\proof Let $A$ be an ample divisor on $\tPP_3$ and let
$\oA=\tpi^\ast(A)$. By the Nakai-Moishezon criterion and
the projection formula, $\oA$ is ample on $\tV$. But all
of $\cohomtV$ is algebraic \cite[p.~120]{clemens},
so for $l\in\ZZ$ we can consider
the linear system $\sL=\left|\alpha+l\oA\right|$. For $l\gg 0$ both
$\sL$ and $\left.\sL\right|_{\oB}$ are free. Hence the general element
$T\in\sL$ is a smooth surface intersecting $\oB$ transversal in a smooth curve
$\oD$. Now fix such a $T$ and let $\oS=T+\ttau(T)$, $\tS=\tpi(\oS)$ and
$\tD=\tpi(\oD)$. But $\left.\tpi\right|_{\oB}$ is an isomorphism, so
$\tD$ is also smooth. Now $\tpi^\ast(\tS).\oB=(T+\ttau(T)).\oB=2\oD$,
hence
\begin{equation*}
    2\tD=2\tpi_\ast(\oD)=
    \tpi_\ast(\tpi^\ast(\tS).\oB)=\tS.\tpi_\ast(\oB)=\tS.\tB
\end{equation*}
by the projection formula. Outside of $\oD$ both maps
$\left.\tpi\right|_T\colon T\rightarrow\tS$ and
$\left.\tpi\right|_{\ttau(T)}\colon\ttau(T)\rightarrow\tS$ are
biholomorphic. Local computations show that these maps are also
biholomorphic on $\oD$. In particular $\tS$ is a smooth surface.
Now let $S=\sigma(\tS)$, $D=\sigma(\tD)$. Since $\tD$ is smooth it cannot
contain any exceptional component, so $\tD$ is the proper transform
of $D$ and $S.B=2D$. Hence $S,D$ and $\tS,\tD$ satisfy A1 and A2 of
the proposition. By construction
$[\oS]=[T+\ttau(T)]=\alpha+[\oA]+\ttau^\ast(\alpha+[\oA])$,
so also A3 is satisfied.
\proofend
Let $\alpha\in C$. We will call any surface $S\subset\PP_3$ which satisfies
A1, A2 and A3 of proposition \ref{proposition:alpha}
an $\alpha$-surface. Let $S$ and
$S'$ be $\alpha$-surfaces. Then $[\oS]=\beta+\ttau^\ast(\beta)$,
$[\oSd]=\beta'+\ttau^\ast(\beta')$, where $\beta=\alpha+\gamma$,
$\beta'=\alpha+\gamma'$, $\gamma, \gamma'\in I$. Hence
$[\oS+\oSd]=[\oS]+[\oSd]=2(\alpha+\ttau^\ast(\alpha)) + 2(\gamma+\gamma')\in 2I$.
This implies $[\tS+\tS']\in 2\cohomtPP$, hence by lemma
\ref{lemma:cut} $S+S'$ cuts out the
empty even set of nodes. So both $S$ and $S'$ cut out the same even
set of nodes $w(\alpha)\in\oC_B$. So we can define a map of sets
\begin{equation*}
    \begin{array}{r@{}ccc}
        w\colon & C & \longrightarrow & \oC_B \\
        & \alpha & \longmapsto & w(\alpha)
    \end{array}
\end{equation*}
which assigns to $\alpha\in C$ the unique even set of nodes
determined by proposition \ref{proposition:alpha}.
The $\alpha$-surfaces are nice surfaces
cutting out even sets of nodes on $B$: They are normal and
their singular locus is contained in $\sS$. Every singularity
of an $\alpha$-surface can be resolved in one blow up.
These singularities are called superisolated.
We will show later that for an $\alpha$-surface $S$ not only
$[\oS]\in\cohomtV$ splits, but also $\oS$ splits into
two smooth surfaces. Another particularly nice class of
contact surfaces are nodal surfaces $T$, but here almost never
$\oT$ splits.
\begin{lemma}
    \begin{itemize}
        \item[i)] The map $w$ is a homomorphism of $\ZZ$-modules.
        \item[ii)] The induced map
            $w_2\colon C\otimes\FF\rightarrow\oC_B$ is an injective
            map of vector spaces.
    \end{itemize}
\end{lemma}\label{lemma:injective}
\proof i) It suffices to show that $w(\alpha_1+\alpha_2)=w(\alpha_1)
+w(\alpha_2)$ for all $\alpha_1,\alpha_2\in C$. So let $\alpha_i\in C$
and let $S_i$ be an $\alpha_i$-surface, $i=1,2$. Moreover let
$S$ be an $(\alpha_1+\alpha_2)$-surface. Then
\begin{alignat*}{2}
    [\oS_i] &=\alpha_i+\ttau^\ast(\alpha_i)+2\gamma_i, & \quad
         & \gamma_i\in I, i=1,2,\\
    [\oS]   &=\alpha_1+\alpha_2+\ttau^\ast(\alpha_1+\alpha_2)
        +2\gamma,& \quad & \gamma\in I.
\end{alignat*}
Hence $[\oS+\oS_1+\oS_2]\in 2I$, so $[\tS+\tS_1+\tS_2]\in 2\cohomtPP$.
By lemma \ref{lemma:cut} $S+S_1+S_2$ cuts out
the empty even set of nodes. So $S$ cuts out the same even set of nodes
as $S_1+S_2$, which is $w(\alpha_1)+w(\alpha_2)$.

ii) Let $w=w(\alpha)\in\oC_B$ and let $S'$ be a surface cutting out $w$
which satisfies A1 and A2 of proposition \ref{proposition:alpha}.
Moreover let $S$ be any $\alpha$-surface. Then $S$ cuts out $w$,
so lemma \ref{lemma:cut} implies $[\tS-\tS']\in 2\cohomtPP$, hence
$[\oS-\oSd]\in 2I$. We have $[\oS]=\beta+\ttau^\ast(\beta)$ with
$\beta\equiv\alpha\bmod I$, so
$[\oSd]=\beta+\gamma+\ttau^\ast(\beta+\gamma)$ for some $\gamma\in I$.
Then $\beta+\gamma\equiv\alpha\bmod I$, so $S'$ is an $\alpha$-surface.
 Since i) holds $w_2$ is well defined. Now assume that
$w=w_2(\alpha_1\bmod 2)=w_2(\alpha_2\bmod 2)$ for some $\alpha_1$, $\alpha_2\in C$.
Let $S$ be an $\alpha_1$-surface. Then $S$ is also an
$\alpha_2$-surface. So $[\oS]=\beta+\ttau^\ast(\beta)$ with
$\beta\equiv\alpha_1\bmod I\equiv\alpha_2\bmod I$, so
$\alpha_1=\alpha_2$ in $C$.
\proofend
Let $\oC_B^d=\im(w_2)\cong\FF^d$. We get the following useful characterisation
of the elements of $\oC_B^d$.
\begin{proposition}\label{proposition:char}
    Let $w\in\oC_B$.
    The following statements are equivalent:
    \begin{itemize}
        \item[1)] $w\in\oC_B^d$.
        \item[2)] There exists a surface cutting out $w$ satisfying
                  A1 and A2 of proposition \ref{proposition:alpha}.
        \item[3)] There exists a surface $S$ cutting out $w$ such that
                  $\tS$ is smooth and $\tS.\tB=2\tD'$, where
                  $\tD'$ is an effective (not necessarily smooth or reduced)
                  divisor on $\tB$.
    \end{itemize}
\end{proposition}
\proof
1)$\Rightarrow$2): $w=w(\alpha)$ for some $\alpha\in C$. Then 
any $\alpha$-surface satisfies 2).
2)$\Rightarrow$3) is obvious.
3)$\Rightarrow$1): We will show that $\oS$ splits into two
smooth surfaces. Consider the double cover
$\left.\tpi\right|_{\oS}\colon\oS\rightarrow\tS$. For every point
$P\in\tD'$ we can find holomorphic coordinates $(x,y,z)$ in a small
open neighbourhood $U_P$ of $P$ such that
\begin{equation*}
    P=(0,0,0),\quad
    \tB\cap U_P=\left\{z=0\right\}\quad\text{and}\quad
    \tS\cap U_P=\left\{h^2(x,y)=z\right\}.
\end{equation*}
Locally the map $\tpi$ is given by $(x,y,w)\mapsto (x,y,w^2)$. Thus
\begin{equation*}
    \oS\cap\tpi^{-1}(U_P)=\left\{h^2(x,y)=w^2\right\}
       =\left\{h(x,y)=w\right\}\cup \left\{h(x,y)=-w\right\}.
\end{equation*}
Both $\left\{h(x,y)=w\right\}$ and $\left\{h(x,y)=-w\right\}$ are smooth
since they are graphs of $\pm h(x,y)$.
For every point $P\not\in\tD'$, there exists a small open
neighbourhood $U_P$ of $P$ such that $\tpi^{-1}(U_P)$ is the
disjoint union of two copies of $U_P$. This defines an atlas of a smooth
manifold $T$ endowed with a map to $\tS$ which is an unbranched double cover.
But $S$ has only isolated singularities, so the group
$\cohom{1}{\smash{\tS},\FF}$ parametrising unbranched double covers
of $\tS$ vanishes \cite[p.~18]{milnor}. Hence $T$ splits into two smooth surfaces
and so $\oS$ does: $\oS=\oS_1+\oS_2$. Now take any class
$\alpha\in C$ with $[\oS_1]\equiv \alpha\bmod I$. Then for
any $\alpha$-surface $S'$ we have $[\oS+\oSd]\in 2I$.
Consequently $[\tS+\tS']\in 2\cohomtPP$, so by lemma \ref{lemma:cut}
$S+S'$ cuts out the empty even set of nodes. Hence
both $S$ and $\oSd$ cut out $w(\alpha)=w$. So $w\in\oC_B^d$.
\proofend
Let $w_1,\ldots,w_d$ be a basis of $\oC_B^d$. Then there exist classes
$\alpha_1,\ldots,\alpha_d\in C$ such that $w(\alpha_i)=w_i$.
Let $S_i$ be an $\alpha_i$-surface, $i=1,\ldots,d$.
The proof of proposition \ref{proposition:alpha} shows that
$\oS_i$ splits into two smooth surfaces
$\oS_i=\oSd_i+\oSdd_i$, $i=1,\ldots,d$.
Then as a consequence of lemma \ref{lemma:injective} we get the
\begin{theorem}\label{theorem:classes}
    $\cohomtV$ is generated by the classes of
    $\tpi^\ast\tH$,
    $\tpi^\ast E_P$, $P\in\sS$ and
    $\oSd_i$, $i=1,\ldots,d$.
\end{theorem}
\section{Applications}
Lemma \ref{lemma:dim} assures that $\dim_{\FF}\oC_B\geq d$, where
$d$ is the defect of the double cover $\pi\colon V\rightarrow\PP_3$
branched along $B$, $b=\deg B$ even. If $B$ is explicitely given,
then $d$ can be computed straightforward by (\ref{equation:defect})
using a computer algebra system. It is however not obvious
how to compute the rank of the 2-torsion in $\cohomstV$ or
equivalently $\dim_{\FF}\oC_B$. However for
$b=4,6$ one can often determine $\dim_{\FF}\oC_B$ and finds
surprisingly often $\dim_{\FF}\oC_B=d$.
We will illustrate this observation by giving a sufficient condition
for this equality. This condition will be applied to the case
of double quartics, i.e.~where $B$ is a nodal quartic surface.

First note that if $w\in\oC_B$ is cut out by a smooth surface $S$
via $D$, then $\tS$ meets every exceptional locus $E_P$, $P\in w$, in a 
line. But $\tB$ meets every exceptional locus
in a smooth conic, so we always have $\tS.\tB=2\tD$ and hence
$w\in\oC_B^d$ by proposition \ref{proposition:char}.
\begin{lemma}\label{lemma:quadric}
    Let $w\in\oC_B$ be cut out by a quadric $S$. Then $w\in\oC_B^d$.
\end{lemma}
\proof If $S$ is smooth we are done.
So we have to study two different cases.

Case I: $S=H_1\cup H_2$ consists of two planes. If $H_1=H_2$ then
$w=0\in\oC_B^d$, so let $H_1\neq H_2$ and $L=H_1\cap H_2$.
Since $b$ is even $L\not\subset B$. So both $H_1$ and
$H_2$ have contact to $B$ and cut out weakly even sets of nodes
$w_1, w_2\in\oC_B$. But $H_1$ and $H_2$ are smooth, so
$w_1, w_2\in\oC_B^d$.
Hence $w=w_1+w_2\in\oC_B^d$.

Case I\!I: $S$ is a quadratic cone with vertex $P$. If $P\not\in\sS$,
then Gallarati's formula \cite[lemma 2.3]{catanese}}
says $\ww -1=b(b-2)$, hence $\ww$ is odd.
This is impossible since $w$ is strictly even. So we have $P\in\sS$.
If $P\in w$ then Gallarati's formula says $\ww -1=b(b-2)$,
contradiction. Hence $P\in\sS\setminus w$. Let $D$ be defined by
$S.B=2D$, then $\mult(D,P)=\tD.E_P$ is even. We have proper transforms 
of the form
\begin{align*}
    \tS &\sim_{lin} 2\tH -2E_P -(\text{other exceptional stuff}),\\
    \tB &\sim_{lin} b\tH -2E_P -(\text{other exceptional stuff}).
\end{align*}
Since $S$ is smooth outside $P$ we also have $\tS.\tB=2\tD +kC_P$ for
some integer $k$. Then $\tS.\tB.E_P=-2C_P.E_P=4$, so
$2k=-kC_P.E_P=2\tD.E_P-4$. It follows that $k$ is even, so
$w\in\oC_B^d$ by proposition \ref{proposition:char}.
\proofend
\begin{proposition}\label{proposition:condition}
    The following statements are equivalent:
    \begin{itemize}
        \item[i)] $\cohomstV$ has no $2$-torsion.
        \item[ii)] $\dim_{\FF}\oC_B =d$.
        \item[iii)] There exists a basis $w_1,\ldots,w_c$ of
            $\oC_B$ such that $w_i$ is cut out by a surface
            $S_i$ which is either smooth, a quadric or
            has the following properties: $\tS_i$ is smooth
            and $\tS_i.\tB=2\tD_i$ for an effective (not necessarily smooth
            or reduced) divisor $\tD_i$ on $\tB$.
    \end{itemize}
\end{proposition}
\proof Combine lemma \ref{lemma:dim} with proposition
\ref{proposition:char} and lemma \ref{lemma:quadric}.\proofend
Now let $B$ be a quartic surface with $\mu$ nodes $\sS=\{P_1,\ldots,P_\mu\}$.
It is a classical result that $\mu\in\{0,\ldots,16\}$ and that all
values really do occur \cite{kummer,rohn}. We will compute the rank of the
$2$-torsion in $\cohomstV$ by computing the code of $B$. First
we give a list of all codes that can appear.

For any even set of nodes $w\in\oC_B$, its {\em weight} $\ww$
denotes the number of nodes it contains.
The elements of any (binary linear) code $C$
will often be called words. There is a natural
partition of $C$ in to the sets $C_d=\{w\in C\mid\ww =d\}$.
We will use two notations for codes in the sequel.
A $[n,k,d]$-code $C$ is a $k$-dimensional subvectorspace of
$\FF^n$ such that any nonzero word $w\in C\setminus\{0\}$ has
weight $\ww\geq d$. A 
$[n,k,\{{d_1}_{m_1},\ldots,{d_l}_{m_l}\}]$-code $C$ is a 
$k$-dimensional subvectorspace of $\FF^n$ with weights
$\{d_1,\ldots,d_l\}=\{\ww\mid w\in C\setminus\{0\}\}$.
Furthermore $C$ is equipped with
a basis $w_1,\ldots,w_k$ such that exactly $m_i$ elements
of the basis have weight $d_i$, $i=1,\ldots,l$. A subscript
$m_i=0$ will be omitted.

Nearly all codes of nodal quartic surfaces can be read off from
Rohn's classification \cite{rohn}. However we will
give a proof which uses only elementary coding theory.
\begin{theorem}\label{theorem:qcodes}
    (Rohn)
    Let $B\subset\PP_3$ be a quartic surface with $\mu$ nodes
    as its only singularities.
    The possibilities for $\oC_B$ being nonzero 
    are given by the following table:
    \begin{equation*}
        \begin{array}{|r||l|}\hline
            \mu & \text{possible codes}\\\hline\hline
              6,7 & [\mu,1,\{6_1\}] \\\hline
              8,9 & [\mu,1,\{6_1\}], [\mu,1,\{8_1\}]\\\hline
            10,11 & [\mu,1,\{10_1\}], [\mu,1,\{8_1\}], [\mu,1,\{6_1\}],
                    [\mu,2,\{6_2,8\}]\\\hline
               12 & [12,2,\{6_1,8_1,10\}], [12,2,\{8_2\}], [12,3,\{6_3,8\}],
                    [12,2,\{6_2,8\}]\\\hline
               13 & [13,3,\{6_3,8\}], [13,3,\{6_3,8,10\}]\\\hline
               14 & [14,4,\{6_4,8,10\}]\\\hline
               15 & [15,5,\{6_5,8,10\}]\\\hline
               16 & [16,6,\{6_6,8,10,16\}]\\\hline
        \end{array}
    \end{equation*}
    Furthermore for every above possibility, $\oC_B$ is unique
    up to a permutation of nodes.
\end{theorem}
\proof
Let $\oC=\oC_B$ and let $C=\{w\in\oC\mid \text{$w$ is strictly even}\}$.
All nonzero weakly even sets of nodes $w\in\oC$ have cardinality
$\ww\in\left\{6,10\right\}$ and all strictly even sets of
nodes $w\in C$ have $\ww\in\left\{8,16\right\}$
\cite[pp.~46--49]{gallarati}. The sum of two strictly
even sets is strictly even,
the sum of two weakly even sets is strictly even and the sum of a 
strictly end a weakly even set is weakly even. Hence $C$ is a subcode
of $\oC$ with 
\begin{equation}\label{align:dim}
\begin{split}
    \dimF C  &\leq\dimF\oC\leq\dimF C+1,\\
    \dimF C  &\geq \mu-11,\\
    \dimF\oC &\geq\mu-10.
\end{split}
\end{equation}
The last two inequalities are due to Beauville \cite[p.~210]{beauville}.
To $\oC$ we can associate a table which we will call weight (addition) table:
\begin{equation*}
    \begin{array}{|c||c|c|c|c|}\hline
               &  6     &    8 &   10 & 16 \\\hline\hline
             6 &  8(10) & {\displaystyle 6(10)\atop\displaystyle 10(12)} &
                          {\displaystyle 8(12)\atop\displaystyle 16(16)} &
              10(16) \\\hline
             8 &        & {\displaystyle 8(12)\atop\displaystyle 16(16)} &
                          {\displaystyle 6(12)\atop\displaystyle 10(14)} & 
               8(16) \\\hline
            10 &        &      &    8(14) &  6(16) \\\hline
    \end{array}
\end{equation*}
All numbers not in brackets denote weights. For example if $w_6\in\oC_6$ and
$w_8\in\oC_8$, then $w_6+w_8$ is weakly even, hence
$w_6+w_8\in\{6,10\}$. In the first case both words overlap in
$4$ nodes, whereas in the second case both words overlap in 2 nodes.
This implies $\mu\geq 10$ in the first case and $\mu\geq 12$ in the second case.
So the numbers in brackets denote the smallest $\mu$ where such an addition could
take place.
Now assume that $\oC\neq\{0\}$. We consider eight different cases.

{\bf Case I:} $\mu\in\{6,7\}$. The weight table exhibits $\dimF\oC =1$.
The only possibility is $\oC=[\mu,1,\{6_1\}]$.

{\bf Case I\!I:} $\mu\in\{8,9\}$. Again  $\dimF\oC=1$, so we
have the two possibilities
$\oC\in\{[\mu,1,\{6_1\}],[\mu,1,\{8_1\}]\}$.

{\bf Case I\!I\!I:} $\mu\in\{10,11\}$. If $\dimF\oC =1$
we have the three possibilities $\oC\in\left\{
[\mu,1,\{6_1\}],[\mu,1,\{8_1\}],[\mu,1,\{10_1\}]\right\}$.
But $C$ is a $[\mu,\dimF C,8]$-code, so the weight table implies
$\dimF C\leq 1$.
Hence  $\dimF\oC\leq 2$. If we have equality then $\oC$
cannot contain a word of weight $10$ by the weight table, hence
$\oC=[\mu,2,\{6_2,8\}]$. This code is up to permutation
of nodes generated by the rows of the following table.
\begin{equation*}
    \newcommand{\rr}{\smash{\hspace*{0.06cm}\blacksquare\hspace*{0.06cm}}}
    \newcommand{\nl}{\\\hline}
    [\mu,2,\{6_2,8\}]:\quad
    \begin{array}{*{10}{|@{}c@{}}|}\hline
        P_1&P_2&P_3&P_4&P_5&P_6&P_7&P_8&P_9&P_{10}\nl\hline
        \rr&\rr&\rr&\rr&\rr&\rr&   &   &   &   \nl
        \rr&\rr&   &   &   &   &\rr&\rr&\rr&\rr\nl
    \end{array}
\end{equation*}

{\bf Case I\!V:} $\mu=12$. Here  $\dimF C\geq 1$ by (\ref{align:dim}) and
$\dim C\leq 2$ by the Griesmer bound
\cite[5.6.2]{lint}. Furthermore $\oC$ contains at
most one word of weight 10 by the weight table.
If $\dim C=1$ we find $\dim\oC=2$ using (\ref{align:dim}).
So $\oC$ contains at least one word of weight 6. If $\oC$ does not contain
a word of weight 10 we find $\oC=[12,2,\{6_2,8\}]$ as before.
If $\oCB$ contains a word of weight $10$ then 
we must have $\oC=\{0,w_6,w_8,w_{10}\}$, $w_i\in\oC_i$, $i=6,8,10$ with
$w_6+w_8=w_{10}$. Hence $\oC$ is (up to permutation) given by the following
table.
\begin{equation*}
    \newcommand{\rr}{\smash{\hspace*{0.06cm}\blacksquare\hspace*{0.06cm}}}
    \newcommand{\nl}{\\\hline}
    [12,2,\{6_1,8_1,10\}]:\quad
    \begin{array}{*{12}{|@{}c@{}}|}\hline
        \rr&\rr&\rr&\rr&\rr&\rr&   &   &   &   &   &   \nl
        \rr&\rr&   &   &   &   &\rr&\rr&\rr&\rr&\rr&\rr\nl
    \end{array}
\end{equation*}
Now let $\dimF C=2$, hence $C=\left<w_8,w'_8\right>$ is generated by
two words $w_8, w_8'$ of weight 8. Up to
permutation, the nonzero words of $C$ are
given by
\begin{equation*}
    \newcommand{\rr}{\smash{\hspace*{0.06cm}\blacksquare\hspace*{0.06cm}}}
    \newcommand{\nl}{\\\hline}
    [12,2,\{8_2\}]:\quad
    \begin{array}{|l|*{12}{|@{}c@{}}|}\hline
        w_8      &\rr&\rr&\rr&\rr&\rr&\rr&\rr&\rr&   &   &   &   \nl
        w'_8     &\rr&\rr&\rr&\rr&   &   &   &   &\rr&\rr&\rr&\rr\nl
        w_8+w'_8 &   &   &   &   &\rr&\rr&\rr&\rr&\rr&\rr&\rr&\rr\nl
    \end{array}
\end{equation*}
Either $\oC=C=[12,2,\{8_2\}]$ or $\dim\oC =3$. Assume that there
exists a $w_{10}\in\oC_{10}$. Let
\begin{align*}
    a_1 &=\left|w_{10}\cap\left\{P_1,\ldots,P_4\right\}\right|, \\
    a_2 &=\left|w_{10}\cap\left\{P_5,\ldots,P_8\right\}\right|, \\
    a_3 &=\left|w_{10}\cap\left\{P_9,\ldots,P_{12}\right\}\right|.
\end{align*}
But $\oC_{10}\leq 1$, so for all $w\in\oC_8$ have
$\left|w+w_{10}\right|=6$, hence $\left|w\cap w_{10}\right|=6$.
Hence $(a_1,a_2,a_3)$ satisfy
\begin{equation*}
    \begin{pmatrix}1&1&0\\ 1&0&1\\ 0&1&1\\ 1&1&1\end{pmatrix}
    \begin{pmatrix}a_1\\a_2\\a_3\end{pmatrix} =
    \begin{pmatrix}6\\6\\6\\10\end{pmatrix}.
\end{equation*}
But this system of linear equations has no solution, so in fact
$\oC_{10}=\emptyset$. Now let $w_6\in\oC_6$. All
$w\in\oC_8$ have
$\left|w+w_6\right|=6$, hence $\left|w\cap w_6\right|=4$.
The corresponding system of equations for $w_6$
has $\!\!{\ }^t\!\left(4,4,4,6\right)$ on the right hand side and we get
$a_1=a_2=a_3=2$. Let $w'_6=w_6+w_8$, $w''_6=w_6+w'_8$. Then
$\oCB=\left<w_6,w'_6,w''_6\right>$ is given (up to permutation) by
\begin{equation*}
    \newcommand{\rr}{\smash{\hspace*{0.06cm}\blacksquare\hspace*{0.06cm}}}
    \newcommand{\nl}{\\\hline}
    [12,3,\{6_3,8\}]:\quad
    \begin{array}{*{12}{|@{}c@{}}|}\hline
        \rr&\rr&\rr&\rr&\rr&\rr&   &   &   &   &   &   \nl
        \rr&\rr&   &   &   &   &\rr&\rr&\rr&\rr&   &   \nl
           &   &\rr&\rr&   &   &\rr&\rr&   &   &\rr&\rr\nl
    \end{array}
\end{equation*}

{\bf Case V:} $\mu=13$. This time we get $\dim C\geq 2$ from
(\ref{align:dim}) and $\dim C\leq 2$ from the Griesmer bound.
Hence $\dim C=2$ and $\dim\oC =3$. As in case V we find that
$|\oC_{10}|\leq 1$ and $|\oC_6|\geq 1$. We can assume that
(after a permutation) $C$ is exactly the code $[12,2,\{8_2\}]$
given in case I\!V. Let us first consider the case $\oC_{10}=\emptyset$.
Let $w_6\in\oC_6$ and define $a_1$, $a_2$ and $a_3$ as in
case V. Let $a_4=\left|w_6\cap\{P_{13}\}\right|$. Then 
the $a_i$'s satisfy
\begin{equation*}
    \begin{pmatrix}1&1&0&0\\ 1&0&1&0\\ 0&1&1&0\\ 1&1&1&1\end{pmatrix}
    \begin{pmatrix}a_1\\a_2\\a_3\\a_4\end{pmatrix}=
    \begin{pmatrix}4\\4\\4\\6\end{pmatrix}.
\end{equation*}
Solving the equations we find $a_1=a_2=a_3=2$ and $a_4=0$, hence
$\oC =[12,3,\{6_3,8\}]$ is just the code from case I\!V.

However if $w_{10}\in\oC_{10}\neq\emptyset$, the same system of
equations for $w_{10}$ has the vector
$\!\!{\ }^t\!\left(6,6,6,10\right)$ on the
right hand side, hence $a_1=a_2=a_3=3$, $a_4=1$. Let
$w_6=w_{10}+w_8$, $w'_6=w_{10}+w'_8$ and $w''_6=w_{10}+w_8+w'_8$. Then
$\oC=\left<w_6,w'_6,w''_6\right>$ is (up to permutation)
given by the following table.
\begin{equation*}
    \newcommand{\rr}{\smash{\hspace*{0.06cm}\blacksquare\hspace*{0.06cm}}}
    \newcommand{\nl}{\\\hline}
    [13,6,\{6_3,8,10\}]:\quad
    \begin{array}{*{13}{|@{}c@{}}|}\hline
        \rr&\rr&\rr&\rr&\rr&\rr&   &   &   &   &   &   &   \nl
        \rr&\rr&   &   &   &   &\rr&\rr&\rr&\rr&   &   &   \nl
        \rr&   &\rr&   &   &   &\rr&   &   &   &\rr&\rr&\rr\nl
    \end{array}
\end{equation*}

{\bf Case VI:}   $\mu=14$. Here we find $\dimF C=3$ and $\dimF\oC =4$.
Furthermore the Griesmer bound implies $\oC_6\neq\emptyset$.
$C$ is generated by three words $w_8, w_8', w_8''$ of weight 8.
Up to permutation of columns, $C$ is given by the following table:
\begin{equation*}
    \newcommand{\rr}{\smash{\hspace*{0.06cm}\blacksquare\hspace*{0.06cm}}}
    \newcommand{\nl}{\\\hline}
    \begin{array}{|l|*{14}{|@{}c@{}}|}\hline
        w_8           &\rr&\rr&\rr&\rr&\rr&\rr&\rr&\rr&   &   &   &   &   &   \nl
        w'_8          &\rr&\rr&\rr&\rr&   &   &   &   &\rr&\rr&\rr&\rr&   &   \nl
        w_8+w'_8      &   &   &   &   &\rr&\rr&\rr&\rr&\rr&\rr&\rr&\rr&   &   \nl
        w_8''         &\rr&\rr&   &   &\rr&\rr&   &   &\rr&\rr&   &   &\rr&\rr\nl
        w_8+w_8''     &   &   &\rr&\rr&   &   &\rr&\rr&\rr&\rr&   &   &\rr&\rr\nl
        w_8'+w_8''    &   &   &\rr&\rr&\rr&\rr&   &   &   &   &\rr&\rr&\rr&\rr\nl
        w_8+w_8'+w_8''&\rr&\rr&   &   &   &   &\rr&\rr&   &   &\rr&\rr&\rr&\rr\nl
    \end{array}
\end{equation*}
Now let $w_6\in\oC_6$. Then for any $w\in\oC_8$ we have
$\left|w_6+w\right|\in\{6,10\}$, hence $w_8$ and $w$ overlap
in 2 or 4 nodes. Let $a_i=\left|w_6\cap\{P_{2i-1},P_{2i}\}\right|$, 
$i=1,\ldots,7$. Then $(a_1,\ldots,a_7)$ satisfy 
\begin{equation*}
    \begin{pmatrix}
        1&1&1&1&0&0&0\\
        1&1&0&0&1&1&0\\
        0&0&1&1&1&1&0\\
        1&0&1&0&1&0&1\\
        0&1&0&1&1&0&1\\
        0&1&1&0&0&1&1\\
        1&0&0&1&0&1&1\\
        1&1&1&1&1&1&1
    \end{pmatrix}
    \begin{pmatrix}
        a_1\\a_2\\\vdots\\a_7
    \end{pmatrix}
    = 
    \begin{pmatrix}
        b_1\\\vdots\\b_7\\6
    \end{pmatrix},
    \quad b_i\in\{2,4\}.
\end{equation*}
Now let $k$ of the $b_i$'s be 4 and $7-k$ of the $b_i$'s be 2.
Then $4\sum a_i =4|w_6|=24=\sum b_i =14+2k$, hence $k=5$.
Thus $|\oC_6|=1+5=6$ and $|\oC_{10}|=2$, and we have
$\left(7\atop 2\right)=21$ cases where this linear system 
of equations has a solution.
In every solution, exactly one of the $a_i$'s has value 2,
exactly four of the $a_i$'s have value 1 and two $a_i$'s vanish.
After a permutation of columns we can assume that
$a_1=2$, $a_2=a_3=1$ and $a_4=0$. Then the second line implies
$a_5+a_6=1$, $a_7=1$.
The fourth line implies $a_5=0$, hence $a_6=1$. So in fact all the
21 solutions are the same up to permutation and $\oC$ is given by
\begin{equation*}
    \newcommand{\rr}{\smash{\hspace*{0.06cm}\blacksquare\hspace*{0.06cm}}}
    \newcommand{\nl}{\\\hline}
    [14,4,\{6_1,8_3,10\}]:\quad
    \begin{array}{*{14}{|@{}c@{}}|}\hline
        \rr&\rr&\rr&\rr&\rr&\rr&\rr&\rr&   &   &   &   &   &   \nl
        \rr&\rr&\rr&\rr&   &   &   &   &\rr&\rr&\rr&\rr&   &   \nl
        \rr&\rr&   &   &\rr&\rr&   &   &\rr&\rr&   &   &\rr&\rr\nl
        \rr&\rr&\rr&   &\rr&   &   &   &   &   &\rr&   &\rr&   \nl
    \end{array}.
\end{equation*}
Let $w_6$ be the word corresponding to the last row of the above
table. Then $\oC=\left<w_8'+w_6,w_6,w_8''+w_6,w_8+w_6\right>$ is given by
the following table.
\begin{equation*}
    \newcommand{\rr}{\smash{\hspace*{0.06cm}\blacksquare\hspace*{0.06cm}}}
    \newcommand{\nl}{\\\hline}
    [14,4,\{6_4,8,10\}]:\quad
    \begin{array}{*{14}{|@{}c@{}}|}\hline
        \rr&\rr&\rr&\rr&\rr&\rr&   &   &   &   &   &   &   &   \nl
        \rr&\rr&   &   &   &   &\rr&\rr&\rr&\rr&   &   &   &   \nl
           &   &\rr&\rr&   &   &\rr&\rr&   &   &\rr&\rr&   &   \nl
        \rr&   &\rr&   &   &   &   &   &\rr&   &\rr&   &\rr&\rr\nl
    \end{array}
\end{equation*}
{\bf Case V\!I\!I}:  $\mu=15$. Here $\dimF C=4$ and $\dimF\oC=5$. Pick one
node $P\in\sS$ and let $\oC'=\{w\in\oC\mid P\not\in w\}$ be the subcode of
$\oC$ of words not containing $P$. For any $w\in\oC\setminus\oC'$ we see
that $\oC$ is spanned by $w$ and $\oC'$. Hence
$\oC$ has dimension 4 and weights in $\{6,8,10\}$.
By case V\!I, $\oC'=[14,4,\{6_4,8,10\}]$.
One can choose $P$ such that there exists a word of weight 6 containing
$P$, so $\oC$ is generated by words of weight 6. It suffices to show
that $\oC$ is unique up to permutation. By similar
(but more tedious) computations as in case V\!I one finds that
$\oC$ is up to permutation given by the following table.
\begin{equation*}
    \newcommand{\rr}{\smash{\hspace*{0.06cm}\blacksquare\hspace*{0.06cm}}}
    \newcommand{\nl}{\\\hline}
    [15,5,\{6_5,8,10\}]:\quad
    \begin{array}{*{15}{|@{}c@{}}|}\hline
        \rr&\rr&\rr&\rr&\rr&\rr&   &   &   &   &   &   &   &   &   \nl
        \rr&\rr&   &   &   &   &\rr&\rr&\rr&\rr&   &   &   &   &   \nl
           &   &\rr&\rr&   &   &\rr&\rr&   &   &\rr&\rr&   &   &   \nl
        \rr&   &   &   &\rr&   &\rr&   &   &   &\rr&   &\rr&\rr&   \nl
           &\rr&\rr&   &   &   &   &   &\rr&   &\rr&   &\rr&   &\rr\nl
    \end{array}
\end{equation*}
{\bf Case V\!I\!I\!I:} $\mu=16$. Here $\dimF C=5$ and $\dimF\oC=6$. Pick one
node $P\in\sS$ and let $\oC'=\{w\in\oC\mod P\not\in w\}$ as before.
Then $\oC'=[14,5,\{6_5,8,10\}]$ by case V\!I\!I.
As before we see that $\oC$ is generated by words of weight 6.
The uniqueness of $\oC$ follows from the fact that it contains a word
$w_{16}\in\oC_{16}$: $\oC$ is spanned by $\oC'$ and $w_{16}$. So we can
generate $\oC$ through words of weight 6 by $\oC'$ and $w_{16}+w_{10}$
for any $w_{10}\in\oC'$. The result is the following table:
\begin{equation*}
    \newcommand{\rr}{\smash{\hspace*{0.06cm}\blacksquare\hspace*{0.06cm}}}
    \newcommand{\nl}{\\\hline}
    [16,6,\{6_6,8,10,16\}]:\quad
    \begin{array}{*{16}{|@{}c@{}}|}\hline
        \rr&\rr&\rr&\rr&\rr&\rr&   &   &   &   &   &   &   &   &   &   \nl
        \rr&\rr&   &   &   &   &\rr&\rr&\rr&\rr&   &   &   &   &   &   \nl
           &   &\rr&\rr&   &   &\rr&\rr&   &   &\rr&\rr&   &   &   &   \nl
        \rr&   &   &   &\rr&   &\rr&   &   &   &\rr&   &\rr&\rr&   &   \nl
           &\rr&\rr&   &   &   &   &   &\rr&   &\rr&   &\rr&   &\rr&   \nl
        \rr&   &\rr&   &   &   &   &   &\rr&   &   &\rr&   &\rr&   &\rr\nl
    \end{array}
\end{equation*}
\proofend
\begin{remark}
    The list of theorem \ref{theorem:qcodes} is compiled such that
    it corresponds directly to the cases denoted
    X\!Ia, X\!Ib, X\!Ic, X\!Id ($\mu=11$),
    X\!I\!Ia, X\!I\!Ib, X\!I\!Ic, X\!I\!Id ($\mu=12$) and
    X\!I\!I\!Ia, X\!I\!I\!Ib ($\mu=13$)
    in Rohn's work.
\end{remark}
Now it is quite easy to compute the rank of the 2-torsion in
$\cohomstV$. Proposition 3.10 in \cite{endrassnodes} tells us that all even
sets of 6 nodes (resp.~8 nodes) on $B$ are cut out by planes
(resp.~reduced quadrics). Hence if $\oC_B$ is
generated by words of weight 6 and 8, then $\cohomstV$ has
no 2-torsion. In the list of theorem \ref{theorem:qcodes}
there are only two cases left: $\mu\in\{10,11\}$ and
$\oC_B=[\mu,1,\{10_1\}]$.
By \cite[theorem~2.23]{catanese} such a quartic $B$ is linearly symmetric,
i.e.~$B=\{\det A=0\}$ where $A$ is a symmetric $4\times 4$
matrix of linear forms. The matrix $A$ is not unique, but every
conjugacy class $[A]$ corresponds to a unique even set of
10 nodes $w_{10}\in\oC_B$ such that (as sets)
$w_{10}$ is exactly the locus where all $3\times 3$ minors
of $A$ vanish. A linearly symmetric quartic
is called symmetroid.
Now by (\ref{equation:defect})
the defect of $V$ is given by
\begin{equation*}
    d = \left\{\begin{array}{l@{\quad}l}
        \dim\sM    & \text{if $\mu=10$,}\\
        \dim\sM +1 & \text{if $\mu=11$.}
    \end{array}\right.
\end{equation*}
Here $\sM$ is the space
of quadratic polynomials vanishing in all nodes of $B$.
If $\mu=11$ then $1\leq d\leq \dimF\oC_B =1$, so
$\sM=\{0\}$ and $\cohomstV$ has no $2$-torsion.
However in the case $\mu=10$ we see that $\cohomstV$ has no $2$-torsion
iff the $10$ nodes of $B$ lie on a quadric surface.
The following lemma we learnt from D.~van Straten.
\begin{lemma}\label{lemma:10}
    Let $B$ be a nodal quartic surface admitting an
    even set of $10$ nodes $w_{10}\in\oC_B$.
    Then the points of $w_{10}$ do not lie
    on a quadric surface.
\end{lemma}
\proof Let $P=\CC[a_1,\ldots,a_{10},x,y,z,w]$ and let
\begin{equation*}
    A = \begin{pmatrix}
        a_1&a_2&a_3&a_4\\
        a_2&a_5&a_6&a_7\\
        a_3&a_6&a_8&a_9\\
        a_4&a_7&a_9&a_{10}
    \end{pmatrix}
\end{equation*}
be the ``general'' symmetric $4\times 4$ matrix.
Define $I=\minors(3,A)$ to be the ideal generated by the
$3\times 3$ minors of $A$. If we cut the variety
$V=V(I)$ with 10 general hyperplane sections
\begin{equation}\label{equation:hs}
    h_i=a_i+\alpha_ix+\beta_iy+\gamma_iz+\delta_iw=0,
    \quad \alpha_i, \beta_i, \gamma_i, \delta_i\in\CC,
    \quad i=1,\ldots,10,
\end{equation}
we obtain the 10 nodes of a symmetroid via the
\newcommand{\hs}{(h_1,\ldots,h_{10})}
identification $P/\hs\cong S=\CC[x,y,z,w]$. So let
$J=\hs$ be generated by 10 hyperplanes as in (\ref{equation:hs})
such that $V(I+J)=\{P_1,\ldots,P_{10}\}\subset\PP_3$.

Let $R=P/I$ and $R'=P/J$.
Denote by $\oI$ (resp.~$\oJ$) the ideal
$I/(J\cap I)\subseteq R'$ (resp.~$J/(I\cap J)\subseteq R$).
We wish to show first that $\hs$ is a regular
sequence in $R$. But $R$ is Cohen-Macaulay \cite[p.~83]{concini}, so
$\hs$ is a regular sequence iff
$\dim R/\oJ=\dim R - 10$. Clearly
$\dim R/\oJ=\dim P/(I+J)=1$. On the other hand $a_1,\ldots,a_{10}$
are general, so $\codim I=3$ \cite[14.4.11]{fulton}.
Thus $\dim R=11$ and
$\hs$ is a regular sequence. Now we obtain a  free
resolution of the $S$-module $P/(I+J)$ as follows. Take
any free resolution of the $P$-module $R$
(forgetting about the grading for a moment)
\begin{diagram}[height=2em,width=2em]
        0 & \rTo & P^{\alpha_k} & \rTo & \ldots & \rTo & P^{\alpha_1}
          & \rTo & P^{\alpha_0} & \rTo & R & \rTo & 0
\end{diagram}
and consider the exact sequence of complexes
\begin{diagram}[height=2em,width=2em]
        0 & \rTo & P^{\alpha_k} & \rTo & \ldots & \rTo & P^{\alpha_1}
          & \rTo & P^{\alpha_0} & \rTo & R & \rTo & 0\\
          &      & \dInto^{\cdot h_1}& &        &      & \dInto^{\cdot h_1}
          &      & \dInto^{\cdot h_1} &  & \dInto^{\cdot h_1} & & \\
        0 & \rTo & P^{\alpha_k} & \rTo & \ldots & \rTo & P^{\alpha_1}
          & \rTo & P^{\alpha_0} & \rTo & R & \rTo & 0 \\
          &      & \dOnto       &      &   &      & \dOnto
          &      & \dOnto       &      & \dOnto & & \\
        0 & \rTo & P_1^{\alpha_k} & \rTo & \ldots & \rTo & P_1^{\alpha_1}
          & \rTo & P_1^{\alpha_0} & \rTo & R/(h_1) & \rTo & 0
\end{diagram}
where $P_1=P/(h_1)$. We obtain a free resolution of 
the $P_1$-module $R/(h_1)$. Repeating this process for $h_2$
we obtain a free resolution of the $P_2=P/(h_1,h_2)$-module
$R/(h_1,h_2)$ an so on. After 10 steps we end up with a free
resolution of the $P_{10}=P/J\cong S$-module $R/\oJ\cong P/(I+J)$.
The Hilbert series of $R/\oJ$ can be obtained from the
gradings of this resolution. But the gradings are the same for all
choices of the $h_i$'s, hence the Hilbert series of $R/\oJ$ is
independent of the special choice of the hyperplanes and so is
the Hilbert series of $\oI$. A Macaulay computation
gives
\begin{equation*}
    h_{\oI}(t) =\sum_{i\geq 0}\dim_{\CC}(\oI_i)\ t^i
             =\frac{t^3(6t^2-15t+10)}{(t-1)^4}
             =10t^3+25t^4+46t^5+O(t^6).
\end{equation*}
The crucial observation now is that $I$ is a radical ideal
and that $\deg V=10$. This can be computed with Macaulay.
It follows immediately that $I+J$ is radical, thus also
$\oI$ is radical.

So any quadric through the 10 nodes defines an nonzero
element in the degree 2 part of $\oI$. From the Hilbert
series we see that $\dim_{\CC}(\oI_2)=0$, hence there
is no such quadric.
\proofend
So as a consequence from proposition \ref{proposition:condition} 
and lemma \ref{lemma:10} we get the
\begin{theorem}\label{theorem:quartics}
    Let $\pi\colon V\rightarrow\PP_3$ be the double cover branched
    along a quartic surface with $\mu$ nodes. Then $\cohomstV$
    contains no 2-torsion except in the case $\mu=10$,
    $\oC_B=[10,1,\{10_1\}]$. Then the 2-torsion in
    $\cohomstV$ has rank 2.
\end{theorem}

Recently, there has been some interest in sextic surface
with 65 nodes.
\begin{example}
    Let $B=\{f=0\}$ be a 65-nodal sextic.
    Such a sextic has been discovered first in \cite{barth}.
    Then \cite{pettersen} there exists a quartic
    hypersurface $X\subset\PP_4$ with 42 ordinary double points, such that
    for one fixed node $P\in X$, the projection from $P$ has exactly
    $B$ as branch locus.  Let $\tX$ be the big resolution of all
    nodes of $X$. Then $\tX$ can be obtained from the double cover
    $\pi\colon V\rightarrow\PP_3$ in the following way: The nodes
    of $B$ arise from the 41 nodes of $X\setminus P$ and from
    24 lines on $X$ through $P$. So $\tX$ is an intermediate resolution
    of $V$, where $24$ nodes are resolved small and 41 nodes
    are resolved big. Hence $\cohomtX$ has rank $1+41+d$.
    From (\ref{align:dim}) we get $\dimF\oC_B\geq 13$.
    But \cite[theorem 8.1 and section 9]{jafferuberman} tells us
    that $d=\dim_{\FF} \oC_B\leq 13$, so
    $d\leq \dimF\oC_B=13$. On the other hand (\ref{equation:defect})
    gives $d=\dim\sM +9$. But $\sM$ is the space of all polynomials
    of degree 5 vanishing in $\sS$, hence contains the four
    (linearly independent) partials of $f$. So
    $d\geq 4+9=13$ and $\cohomdtX =55$. This rank can be
    computed in a different way
    \cite[ch.~I\nx I]{werner}: $X$ is a 42-nodal hyperquartic,
    so $\cohomdtX =1+d'+42$, where $d'$ can be computed by means of
    \begin{equation*}
        d'=\dim\sM'-\left(\left(7\atop 4\right)-42\right).
    \end{equation*}
    Here $\sM'$ is the space of all homogeneous polynomials of
    degree 3 passing through $\sing X$. We see that
    \begin{itemize}
        \item Through the nodes of any 65-nodal sextic
              $B=\left\{f=0\right\}$
              there passes exactly one three parameter
              family of quintic surfaces. This family is just the
              linear system spanned by the partial derivatives of $f$.
        \item Through the nodes of $X=\left\{g=0\right\}$
              there passes a 4 parameter
              family of hypercubics. Again it is the linear
              system spanned by the partials of $g$.
    \end{itemize}
    In particular, all 65-nodal sextics $B$ have $\dim_{\FF}\oC_B=13$,
    the double solid $V$ branched over $B$ has defect 13 and
    $\cohomstV$ has no 2-torsion.
\end{example}
Now let $B=\left\{f=0\right\}$ be of even degree $b$ and
consider again Beauville's bound
\begin{equation}\label{beauville:formel}
    \dim_{\FF} \oC_B\geq\mu -b_2\left(\smash{\tB}\right)\nx/2 +1
    = \mu -b\left(b^2-4b+6\right)\nx/2.
\end{equation}
The Miyaoka inequality \cite{miyaoka} implies
$\mu\leq (4/9)b\left(b-1\right)^2$, so the right hand side
of (\ref{beauville:formel}) is always negative for $b\geq 18$.
For small values of $b$ however surfaces with
many nodes must have a nontrivial code. The bound is sharp for
$\left(b,s\right)\in\left\{\left(4,16\right),\left(6,65\right)\right\}$.
However this bound is not sharp for $b\geq 8$, as is suggested by
the following example.
\begin{example}\label{example:x8}
    Let $B=X_8$ be the 168-nodal octic of \cite{endrassoctic}.
    A Macaulay computation gives $\dim_{\FF} \oC_B \geq d=19$,
    whereas the bound exhibits $\dim_{\FF} \oC_B\geq 168-151+1=18$.
\end{example}
Now (\ref{beauville:formel}) can be improved using (\ref{equation:defect})
as follows.
The jacobian ideal $J$ of $\left(f\right)$ has a homogeneous
decomposition $J=\bigoplus_{n\geq b-1}J_n$. Clearly
$J_{3b/2-4}\subseteq \sM$, and the expected dimension is
$\dim_{\CC}\left(J_{3b/2-4}\right)=4\left(b/2\atop 3\right)$.
This is however only true if the partial derivatives
$f_x$, $f_y$, $f_z$ and $f_w$ do not satisfy a nontrivial
relation in $J_{3b/2-4}$. For trivial reasons this holds
for $b=2,4,6$. The first obvious relations (in $J_{2(b-1)}$)
$f_x f_y=f_y f_x$ etc.~are called Koszul relations.
\begin{lemma}
    There are no nontrivial relations between
    $f_x$, $f_y$, $f_z$ and $f_w$ in degree $\leq 3b/2-2$.
\end{lemma}
\proof
Suppose 
$\alpha f_x+\beta f_y+\gamma f_z+\delta f_w=0$ is a nontrivial relation
where $\alpha$, $\beta$, $\gamma$, $\delta$ are homogeneous
of degree $k\leq \left(b-2\right)/2$. Consider the linear system
$\PP\left(J_{b-1}\right)$. Its base locus is just $\sS$. Since
the hessian of $f$ in $P\in\sS$ has rank three, three general
elements of $\PP\left(J_{b-1}\right)$ intersect transversal in $P$.
Hence we can assume that (after a projective coordinate change)
the intersection
$Z=\left\{f_x=f_y=f_w=0\right\}$ consists of $\left(b-1\right)^3$
different points.

But the restriction $\left.\delta f_w\right|_Z=0$
vanishes identically. By Miyaoka, $f_w$ vanishes on
at most $(4/9)b\left(b-1\right)^2$ points of $Z$. 
Let $G$, $H$ be general (smooth) elements in the linear system
associated to $f_x$, $f_y$ and $f_z$. Since
$\left\{\delta =0\right\}.G.H=k\left(b-1\right)^2$,
$\delta$ vanishes on at most $\left(b-2\right)\left(b-1\right)^2\nx/2$
points of $Z$. Thus
\begin{equation*}
    9\left(b-2\right)\left(b-1\right)^2+
    8b\left(b-1\right)^2 \geq 18\left(b-1\right)^3,
\end{equation*}
and consequently $b=0$. This proves the lemma.\proofend
So our formula $\dim_{\CC}\left(J_{3b/2-4}\right)=4\left(b/2\atop 3\right)$
is correct, hence we have derived a now bound by replacing
$\dim\sM$ by $\dim J_{3b/2-4}$ in equation (\ref{equation:defect}).
\begin{proposition}\label{proposition:dim}
    Let $B\subset\PP_3$ be a surface of even degree $b$ with
    $\mu$ nodes as its only singularities. Then
    \begin{equation*}
        \dim_{\FF}\oC_B\geq \mu -\left(b-2\right)
        \left(23b^2-38b+24\right)\nx /48.
    \end{equation*}
\end{proposition}
This bound slightly improves Beauville's bound.
For the surface $X_8$ of example \ref{example:x8}, our bound gives
$\dim\oC_{X_8}\geq 19$. 
We do not know an example of a surface of degree $\geq 8$ 
where our bound is sharp.
Also for $b\geq 24$ the right hand side of our bound
is always negative.
Surprisingly, the coefficient $23/48$ of $b^3$ in this bound
is the same as in the spectral upper bound for nodes
\cite[pp.~417-418]{arnold}.
\noindent%
Johannes Gutenberg-Universit\"at\\
Fachbereich 17\\
Staudinger-Weg 9\\
55099 Mainz, Germany\\
{\tt endrass@mathematik.uni-mainz.de}
\end{document}